\documentclass [a4paper,11pt]{article}
\def\dateref{18.06.2008}

\usepackage[latin1]{inputenc}
\usepackage{amsmath}
\usepackage{amssymb}
\usepackage{amsthm}
\usepackage[dvips]{graphicx}

\newtheorem{df}{Definition}
\newtheorem{lm}{Lemma}
\newtheorem{cl}{Corollary}
\newtheorem{theo}{Theorem}
\newtheorem{prop}{Proposition}

\newcommand{\ol}[1] {\overline{#1}}

\def\be{\begin{equation}}
\def\ee{\end{equation}}
\def\bea{\begin{eqnarray}}
\def\eea{\end{eqnarray}}
\def\R {\mathbb{R}} \def\N {\mathbb{N}} \def\Z {\mathbb{Z}} \def\Zp {\Z_+}
\def\diam {\mathrm{diam}\, }
\def\eps{\varepsilon}
\def\phi {\varphi}

\def\Pr {\textbf{Proof: }}
\def\bk {\backslash}
\def\htop {h_{\mathrm{top}}}

 \def\cB {\mathcal{B}} \def\cC {\mathcal{C}} \def\cD {\mathcal{D}}
\def\cE {\mathcal{E}}  \def\cG {\mathcal{G}} 
   
 \def\cN {\mathcal{N}}  
 \def\cR {\mathcal{R}}  
\def\cU {\mathcal{U}}

    \def\Sk {\Sigma_k} 

\def\uu {\underline{u}} \def\uv {\underline{v}}
\def\ux {\underline{x}}  \def\uw {\underline{w}} \def\uz {\underline{z}}
\def\ub {\underline{\beta}} \def\ob {\overline{\beta}}
\def\Tab {T_{\alpha,\beta}}
\def\Sab {\Sigma_{\alpha,\beta}} \def\Suv {\Sigma_{\uu,\uv}}
 \def\muab {\mu_{\ab}} \def\hmuab {\hat{\mu}_{\ab}}
\def\phiab {\phi^{\alpha,\beta}}

\def\uab {\underline{u}^{\alpha,\beta}} \def\vab {\underline{v}^{\alpha,\beta}}
 
\def\uabb {\uu^{\alpha(\beta),\beta}} 

\def\ab {\alpha,\beta} \def\abb {\alpha(\beta),\beta}
\def\A {{\tt A}} \def\As {\A^*}
\def\i {{\tt i}} \def\iab {\i^{\ab}} \def\Tabb {T_{\abb}} \def\muabb {\mu_{\abb}}
\def\phiabb {\phi^{\abb}} 
\def\ib {\i^\beta} \def\Tb {T_\beta} \def\mub {\mu_\beta} \def\hmub {\hat{\mu}_\beta}
\def\phib {\phi^\beta} \def\eb {\underline{\eta}^\beta} \def\Sb {\Sigma_\beta}
\def\ue {\underline{\eta}}

\begin{document}

\title{A point is normal for almost all maps $\beta x + \alpha \mod 1$ or generalized $\beta$-maps.}
\author{B. Faller\footnote{e-mail: bastien.faller@a3.epfl.ch} and
C.-E. Pfister\footnote{e-mail: charles.pfister@epfl.ch}\\
EPF-L, Institut d'analyse et calcul scientifique\\
CH-1015 Lausanne, Switzerland}

\date{\dateref}
\maketitle

\begin{abstract}
We consider the  map $T_{\alpha,\beta}(x):=
\beta x + \alpha \mod 1$, which admits a unique probability measure of maximal
entropy $\muab$. For $x \in [0,1]$, we show that the orbit of $x$
is $\mu_{\alpha,\beta}$-normal for almost all
$(\alpha,\beta)\in[0,1)\times(1,\infty)$ (Lebesgue measure).
Nevertheless we construct analytic curves in
$[0,1)\times(1,\infty)$ along them the orbit of $x=0$ is at most
at one point $\mu_{\alpha,\beta}$-normal.  These curves are
disjoint and they fill the set $[0,1)\times(1,\infty)$. We also
study the generalized $\beta$-maps (in particular the tent map).
We show that the critical orbit $x=1$ is normal with respect to
the measure of maximal entropy for almost all $\beta$.
\end{abstract}

\newpage

\section{Introduction}

In this paper, we consider a dynamical system $(X,d,T)$ where $(X,d)$ is a compact metric
space endowed with its Borel $\sigma$-algebra $\cB$ and $T: X \to X$ is a measurable
application. Let $C(X)$ denote the set of all continuous functions from $X$ into $\R$. The
set $M(X)$ of all Borel probability measures is equipped with the weak$^*$-topology.
$M(X,T) \subset M(X)$ is the subset of all $T$-invariant probability measures. For $\mu \in
M(X,T)$, let $h(\mu)$ denote the measure-theoretic entropy of $\mu$. For all $x\in X$ and
$n\geq1$, the empirical measure of order $n$ at $x$ is \be \nonumber \cE_n(x) :=
\frac{1}{n} \sum_{i=0}^{n-1} \delta_x \circ T^{-i} \in M(X), \ee where $\delta_x$ is the
Dirac mass at $x$. Let $V_T(x) \subset M(X,T)$ denote the set of all cluster points of
$\{\cE_n(x)\}_{n\geq1}$ in the weak$^*$-topology.
\begin{df} \label{df normality}
Let $\mu \in M(X,T)$ be an ergodic measure and $x \in X$. The
orbit of $x$ under $T$ is $\mu$-\textbf{normal}, if $V_T(x) =
\{\mu\}$, ie for all continuous $f \in C(X)$, we have \be
\nonumber \lim_{n\to\infty} \frac{1}{n} \sum_{i=0}^{n-1} f(T^ix) =
\int f d\mu. \ee
\end{df}
By the Birkhoff Ergodic Theorem, $\mu$-almost all points are
$\mu$-normal, however it is difficult to identify a $\mu$-normal
point. This paper is devoted to the study of the normality of
orbits for piecewise monotone continuous applications of the
interval. We consider a family $\{T_\kappa\}_{\kappa \in K}$ of
piecewise monotone continuous applications parameterized by a
parameter $\kappa \in K$, such that for all $\kappa \in K$ there
is a unique measure of maximal entropy $\mu_\kappa$. In our case
$K$ is a subset of $\R$ or $\R^2$. For a given $x\in X$, we
estimate the Lebesgue measure of the subset of $K$ such that the
orbit of $x$ under $T_\kappa$ is $\mu_\kappa$-normal.

For example, let $\Tab: [0,1] \to [0,1]$ be the piecewise monotone continuous application
defined by $\Tab(x) = \beta x + \alpha \mod 1$; here $\kappa = (\ab) \in [0,1) \times
(1,\infty)$. In \cite{Pa64}, Parry constructed a $\Tab$-invariant probability measure
$\muab$ absolutely continuous with respect to Lebesgue measure, which is the unique measure
of maximal entropy. The main result of section \ref{sec beta x + alpha alpha const} is
Theorem \ref{theo selfnormality for Tab}, which shows that for all $x\in[0,1]$ the set \be
\nonumber \cN(x) := \{(\ab) \in [0,1) \times (1,\infty): \text{the orbit of $x$ under
$\Tab$ is $\muab$-normal} \} \ee has full $\lambda^2$-measure, where $\lambda^d$ is the
$d$-dimensional Lebesgue measure. This is a generalization of a theorem of Schmeling in
\cite{Sc97}, where the case $\alpha=0$ and $x=1$ is studied. For the $\beta$-maps, the
orbit of $1$ plays a particular role, so the restriction to $x=1$ considered by Schmeling
is natural. Similarly for $\Tab$, the orbits of $0$ and $1$ are very important. In Theorem
\ref{theo selfnormality for Tab along curves}, we show that there exist curves in the plane
$(\ab)$ defined by $\alpha = \alpha(\beta)$ along which the orbits of $0$ or $1$ are never
$\muab$-normal. The curve $\alpha=0$ is a trivial example of such a curve for the fixed
point $x=0$. In section \ref{sec generalized beta-maps}, we study the generalized
$\beta$-maps introduced by G\'ora \cite{Go07}. A generalized $\beta$-map is similar to a
$\beta$-map, but each lap is replaced by an increasing or decreasing lap of constant slope
$\beta$ according to a sequence of signs. For a given class of generalized $\beta$-maps,
there exists $\beta_0$ such that for all $\beta>\beta_0$, there is unique measure of
maximal entropy $\mu_\beta$ and the set \be \nonumber \{\beta>\beta_0: \text{the orbit of
$1$ under $T_\beta$ is $\mu_\beta$-normal} \} \ee has full $\lambda^1$-measure. Since the
tent maps are generalized $\beta$-maps, we obtain an alternative proof of results of Bruin
in \cite{Br98}.

\section{Preliminaries} \label{sec preliminaries}

Let us define properly the coding for a piecewise monotone continuous application of the
interval. The classical papers are \cite{Re57}, \cite{Pa64} and \cite{Ho79}. We consider
the piecewise monotone continuous applications of the following type. Let $k\geq2$ and $0 =
a_0 < a_1 < \cdots < a_k=1$. We set $\A := \{0,\dots,k-1\}$, $I_0 = [a_0,a_1$), $I_j =
(a_j,a_{j+1})$ for all $j \in 1,\dots,k-2$, $I_{k-1} = (a_{k-1},a_k]$ and $S_0 = \{a_j:
j\in 1,\dots,k-1\}$. For all $j \in \A$, let $f_j: I_j \to [0,1]$ be a strictly monotone
continuous map. A piecewise monotone continuous application $T: [0,1]\bk S_0 \to [0,1]$ is
defined by \be \nonumber T(x) = f_j(x) \quad \text{if } x \in I_j. \ee We will state later
in each specific case how to define $T$ on $S_0$. We set $X_0 = [0,1]$ and for all $n\geq1$
\be \label{eq X_n S_n S} X_n = X_{n-1} \bk S_{n-1} \quad \text{and} \quad S_n = \{x \in
X_n: T^n(x) \in S_0\}, \ee so that $T^n$ is well defined on $X_n$. Finally we define $S =
\bigcup_{n\geq0} S_n$ such that $T^n(x)$ is well defined for all $x \in [0,1]\bk S$ and all
$n\geq0$.

Let $\A$ be endowed with the discrete topology and $\Sk = \A^{\Zp}$ be the product space.
The elements of $\Sk$ are denoted by $\ux = x_0x_1\dots$. A finite string $\uw = w_0 \dots
w_{n-1}$ with $w_j \in \A$ is a \textbf{word}. The \textbf{length} of $\uw$ is $|\uw|=n$.
There is a single word of length 0, the \textbf{empty word} $\eps$. The set of all words is
$\As$. For two words $\uw,\uz$, we write $\uw\ \uz$ for the concatenation of the two words.
For $\ux \in \Sk$, let $\ux_{[i,j)} = x_i\dots x_{j-1}$ denote the word formed by the
coordinates $i$ to $j-1$ of $\ux$. For a word $\uw \in \As$ of length $n$, the
\textbf{cylinder} $[\uw]$ is the set \be \nonumber [\uw] := \{\ux \in \Sk: \ux_{[0,n)} =
\uw\}. \ee The family $\{\,[\uw]: \uw \in \As\}$ is a base for the topology and a
semi-algebra generating the Borel $\sigma$-algebra. For all $\beta>1$, there exists a
metric $d_\beta$ compatible with the topology defined by \be \nonumber d_\beta(\ux,\ux') :=
\begin{cases}
\begin{array}{ll} 0 & \text{if } \ux = \ux' \\ \beta^{-\min\{n\geq0: \ \ux_n \neq \ux'_n\}}
& \text{otherwise.} \end{array} \end{cases} \ee The left shift map $\sigma: \Sk \to \Sk$ is
defined by \be \nonumber \sigma(\ux) = x_1x_2\dots. \ee It is a continuous map. We define a
total order on $\Sk$ denoted by $\prec$. We set \be \nonumber \delta(j) =
\begin{cases} +1 & \text{if $f_j$ is increasing} \\ -1 & \text{if $f_j$ is decreasing}
\end{cases} \ee and for word $\uw$ \be \nonumber \delta(\uw) = \begin{cases} 1 & \text{if
$\uw = \eps$} \\ \delta(w_0) \cdots \delta(w_{n-1}) & \text{otherwise.} \end{cases} \ee Let
$\ux \neq \ux' \in \Sk$ and define $n = \min\{j\geq0: x_j \neq x'_j\}$, then \be \nonumber
\ux \prec \ux' \Leftrightarrow \begin{cases} x_n < x_n' & \text{if } \delta(\ux_{[0,n)}) =
+1 \\ x_n > x_n' & \text{if } \delta(\ux_{[0,n)}) = -1. \end{cases} \ee When all maps $f_j$
are increasing, this is the lexicographic order.

We define the coding map $\i:[0,1]\bk S \to \Sk$ by \be \nonumber \i(x) := \i_0(x) \i_1(x)
\dots \quad \text{with } \i_n(x) = j \Leftrightarrow T^n(x) \in I_j. \ee The coding map
$\i$ is left undefined on $S$. Henceforth we suppose that $T$ is such that $\i$ is
injective. A sufficient condition for the injectivity of the coding is the existence of
$\lambda>1$ such that $|f_j'(x)| \geq \lambda$ for all $x \in I_j$ and all $j \in \A$, see
\cite{Pa64}. This condition is satisfied in all cases considered in the paper. The coding
map is order preserving, ie for all $x,x' \in [0,1]\bk S$ \be \nonumber x < x' \Rightarrow
\i(x) \prec \i(x'). \ee Define $\Sigma_T := \ol{\i([0,1]\bk S)}$. We introduce now the
$\phi$-expansion as defined by Parry. For all $j \in \A$, let $\phi^j: [j,j+1] \to
[a_j,a_{j+1}]$ be the unique monotone extension of $f_j^{-1}: (c,d) \to (a_j,a_{j+1})$
where $(c,d) := f_j\big((a_j,a_{j+1})\big)$. The map $\phi: \Sigma_k \to [0,1]$ is defined
by \be \nonumber \phi(\ux) = \lim_{n\to\infty} \phi^{x_0}\Big(x_0 + \phi^{x_1}\big(x_1 +
\dots + \phi^{x_n}(x_n)\big)\Big). \ee Parry proved that this limit exists if $\i$ is
injective. The map $\phi$ is order preserving. Moreover $\phi|_{\i([0,1]\bk S)} = \i^{-1}$
and for all $n\geq0$ and all $x \in [0,1]\bk S$ \be \label{eq T^n = phi sigma^n i} T^n(x) =
\phi \circ \sigma^n \circ \i(x). \ee If the coding map is injective, one can show that the
map $\phi$ is continuous (see Theorem 2.3 in \cite{FaPf08}). Using the continuity and the
monotonicity of $\phi$, we have $\phi(\Sigma_T) = [0,1]$. Remark that there is in general
no extension of $\i$ on $[0,1]$ such that equation \eqref{eq T^n = phi sigma^n i} is valid
on $[0,1]$. For all $j\in\A$, define \be \nonumber \uu^j := \lim_{x\downarrow a_j} \i(x)
\quad \text{and} \quad \uv^j := \lim_{x\uparrow a_{j+1}} \i(x) \qquad \text{with } x \in
[0,1]\bk S. \ee The strings $\uu^j$ and $\uv^j$ are called critical orbits and (see for
instance \cite{Ho79}) \be \label{eq Sigma_T}\Sigma_T = \{\ux \in \Sk: \uu^{x_n} \preceq
\sigma^n \ux \preceq \uv^{x_n} \ \forall n\geq0\}. \ee Moreover the critical orbits $\uu^j,
\uv^j$ satisfy for all $j\in\A$ \be \label{eq uuj uvj}
\begin{cases} \uu^{u^j_n} \preceq \sigma^n \uu^j \preceq \uv^{u^j_n}
\\ \uu^{v^j_n} \preceq \sigma^n \uv^j \preceq \uv^{v^j_n}
\end{cases} \quad \forall n\geq0. \ee

Let us recall the construction of the Hausdorff dimension. Let $(X,d)$ be a metric space
and $E \subset X$. Let $\cD_\eps(E)$ be the set of all finite or countable cover of $E$
with sets of diameter smaller then $\eps$. For all $s\geq0$, define \be\nonumber
H_\eps(E,s) := \inf\{\sum_{B \in \cC} (\diam B)^s: \cC \in \cD_\eps(E) \} \ee and the
$s$-Hausdorff measure of $E$, $H(E,s) := \lim_{\eps \to 0} H_\eps(E,s)$. The Hausdorff
dimension of $E$ is \be \nonumber \dim_HE := \inf\{s\geq0: H(E,s)=0\}. \ee In \cite{Bo73},
Bowen introduced a definition of the topological entropy of non compact set for a
continuous dynamical system on a metric space. We recall this definition. Let $(X,d)$ be a
metric space, $T: X \to X$ a continuous application. For $n\geq1$, $\eps>0$ and $x \in X$,
let \be \nonumber B_n(x,\eps) = \{y \in X: d(T^j(x),T^j(y)) < \eps \ \forall
j=0,\dots,n-1\}. \ee For $E \subset X$, such that $T(E) \subset E$, let $\cG_n(E,\eps)$ be
the set of all finite or countable covers of $E$ with Bowen's balls $B_m(x,\eps)$ for
$m\geq n$. For all $s\geq0$, define \be \nonumber C_n(E,\eps,s) := \inf\{\sum_{B_m(x,\eps)
\in \cC} e^{-ms}: \cC \in \cG_n(x,\eps)\} \ee and $C(E,\eps,s) := \lim_{n\to\infty}
C_n(E,\eps,s)$. Now, let \be \nonumber \htop(E,\eps) := \inf\{s\geq0: C(E,\eps,s) = 0\} \ee
and finally $\htop(E) = \lim_{\eps\to0} \htop(E,\eps)$ (this last limit increase to
$\htop(E)$). There is an evident similarity of this definition with the Hausdorff
dimension. This similarity is the key of the next lemma.

\begin{lm} \label{lm htop dimH}
For $\beta>1$, consider the dynamical system $(\Sk,d_\beta,\sigma)$. Let $E \subset \Sk$ be
such that $\sigma(E) \subset E$, then \be \nonumber \dim_HE \leq
\frac{\htop(E)}{\log\beta}. \ee
\end{lm}

\Pr Let $\eps\in(0,1), s\geq0, n\geq0$ and $\cC \in \cG_n(E,\eps)$. Since $\diam
B_m(x,\eps) \leq \eps \beta^{-m+1} \leq \eps \beta^{-n+1}$ for all $B_m(x,\eps) \in \cC$,
$\cC$ is a cover of $E$ with sets of diameter smaller than $\eps \beta^{-n+1}$. Moreover
\be \nonumber \sum_{B_m(x,\eps) \in \cC} \diam(B_m(x,\eps))^\frac{s}{\log\beta} \leq
(\eps\beta)^\frac{s}{\log\beta} \sum_{B_m(x,\eps) \in \cC} e^{-ms}. \ee Thus
$H_\delta(E,\frac{s}{\log\beta}) \leq (\eps\beta)^\frac{s}{\log\beta} C_n(E,\eps,s)$ with
$\delta = \eps \beta^{-n+1}$. Taking the limit $n\to\infty$, we obtain \be \nonumber
H(E,\frac{s}{\log\beta}) \leq (\eps\beta)^\frac{s}{\log\beta} C(E,\eps,s). \ee If $s >
\htop(E,\eps)$, then $H(E,\frac{s}{\log\beta}) = 0$ and $\frac{s}{\log\beta} \geq \dim_HE$.
This is true for all $s > \htop(E,\eps)$, thus \be \nonumber \dim_HE \leq
\frac{\htop(E,\eps)}{\log\beta} \leq \frac{\htop(E)}{\log\beta}. \quad \Box \ee

The next lemma is a classical result about the Hausdorff dimension, it is Proposition 2.3
in \cite{Fa03}.

\begin{lm} \label{lm dimH Holder continuity}
Let $(X,d),(X',d')$ be two metric spaces and $\rho: X \to X'$ be
an $\alpha$-H\"older continuous application with $\alpha \in
(0,1]$. Let $E \in X$, then \be \nonumber \dim_H\rho(E) \leq
\frac{\dim_HE}{\alpha}. \ee
\end{lm}

Finally we report Theorem 4.1 from \cite{PfSu07}. This theorem is used to estimate the
topological entropy of sets we are interested in.

\begin{theo} \label{theo htop ^K G}
Let $(X,d,T)$ be a continuous dynamical system and $F \subset M(X,T)$ be a closed subset.
Define \be \nonumber G := \{x \in X: V_T(x) \cap F \neq \emptyset \}. \ee Then \be
\nonumber \htop(G) \leq \sup_{\nu \in F} h(\nu). \ee
\end{theo}

\section{Normality for the maps $\beta x + \alpha \mod 1$}
\label{sec beta x + alpha alpha const}

In this section, we study the piecewise monotone continuous applications $\Tab$ defined by
$\Tab(x) = \beta x + \alpha \mod 1$ with $\beta>1$ and $\alpha \in [0,1)$. These maps were
studied by Parry in \cite{Pa64} as a generalization of the $\beta$-maps. In his paper Parry
constructed a $\Tab$-invariant probability measure $\muab$, which is absolutely continuous
with respect to Lebesgue measure. Its density is \be \label{eq density muab} h_{\ab}(x) :=
\frac{d\muab}{d\lambda} (x) = \frac{1}{N_{\ab}} \ \frac{\sum_{n\geq0} 1_{x<\Tab^n(1)} -
\sum_{n\geq0} 1_{x<\Tab^n(0)}} {\beta^{n+1}}, \ee with $N_{\ab}$ the normalization factor.
In \cite{Ha75}, Halfin proved that $h_{\ab}(x)$ is nonnegative for all $x \in [0,1]$. Set
$k := \lceil \alpha + \beta \rceil$ and let $\iab$ denote the coding map under $\Tab$,
$\phiab$ the corresponding $\phi$-expansion, $\Sab := \Sigma_{\Tab}$, $\uab :=
\lim_{x\downarrow0} \iab(x)$ and $\vab := \lim_{x\uparrow1} \iab(x)$. We specify how $\Tab$
is defined at the discontinuity points. We choose to define $\Tab$ by right-continuity at
$a_j \in S_0$. Doing this we can also extend the definition of the coding map $\iab$ using
the disjoint intervals $[a_j,a_{j+1})$ for $j\in \A$, so that $\iab$ is now defined for all
$x \in [0,1)$ \footnote{This convention differs from that  made in the previous section;
however it is the most convenient choice when all $f_j$ are increasing.}. We can show that
$\uab = \iab(0)$ and \be \nonumber \i([0,1)) = \{\ux \in \Sk: \uab \preceq \sigma^n \ux
\prec \vab \quad \forall n\geq0\} \ee and equation \eqref{eq T^n = phi sigma^n i} is true
for all $x \in [0,1)$. It is easy to check that formula \eqref{eq Sigma_T} becomes \be
\label{eq Sab} \Sab = \{\ux \in \Sk: \uab \preceq \sigma^n \ux \preceq \vab \quad \forall
n\geq0\} \ee and inequations \eqref{eq uuj uvj} become \be \label{eq uab vab}
\begin{cases} \uab \preceq \sigma^n \uab \preceq \vab \\ \uab \preceq \sigma^n \vab \preceq
\vab \end{cases} \quad \forall n\geq0. \ee It is known that the
dynamical system $(\Sab,\sigma)$ has topological entropy
$\log\beta$. Moreover, Hofbauer showed in \cite{Ho80} that it has
a unique measure of maximal entropy $\hat{\mu}_{\ab}$, $\muab =
\hat{\mu}_{\ab} \circ (\phiab)^{-1}$ and $\muab$ is the unique
measure of maximal entropy for $\Tab$. In view of \eqref{eq Sab}
and \eqref{eq uab vab}, for a couple $(\uu,\uv) \in \Sk^2$
satisfying \be \label{eq uu uv}
\begin{cases} \uu \preceq \sigma^n \uu \preceq \uv \\ \uu \preceq \sigma^n \uv \preceq \uv
\end{cases} \quad \forall n\geq0, \ee we define the shift space \be \label{eq Suv} \Suv :=
\{\ux \in \Sk: \uu \preceq \sigma^n \ux \preceq \uv \quad \forall n\geq0\}. \ee

We give now a lemma and a proposition which are the keys of the
main theorem of this section. In the lemma, we show that for given
$x$ and $\alpha$, there is exponential separation between the
orbits of $x$ under the two different dynamical systems
$T_{\ab_1}$ and $T_{\ab_2}$. The proposition asserts that the
topological entropy of $\Suv$ depends continuously on the the
critical orbits $\uu$ and $\uv$.

\begin{lm} \label{lm separation of coding in Tab}
Let $x \in [0,1)$, $\alpha \in [0,1)$ and $1 < \beta_1 \leq \beta_2$. Define $l =
\min\{n\geq0: \i^1_n(x) \neq \i^2_n(x) \}$ with $\i^j(x) = \i^{\ab_j}$ for $j = 1,2$. If
$x\neq0$, then \be \nonumber \beta_2-\beta_1 \leq \frac{\beta_2^{-l+1}}{x}. \ee If $x=0$
and $\alpha\neq0$, then \be \nonumber \beta_2-\beta_1 \leq \frac{\beta_2^{-l+2}}{\alpha}.
\ee
\end{lm}

\Pr Let $\delta := \beta_2-\beta_1 \geq0$. We prove by induction that for all $m\geq1$,
$\i^1_{[0,m)}(x) = \i^2_{[0,m)}(x)$ implies \be \nonumber T_2^m(x) - T_1^m(x) \geq
\beta_2^{m-1}\delta x, \ee where $T_i = T_{\ab_i}$. For $m=1$, \be \nonumber T_2(x)-T_1(x)
= \beta_2 x + \alpha -\i^2_0(x) - (\beta_1 x + \alpha -\i^1_0(x)) = \delta x. \ee Suppose
that this is true for $m$, then $\i^1_{[0,m+1)} = \i^2_{[0,m+1)}$ implies \bea \nonumber
T_2^{m+1}(x) - T_1^{m+1}(x) &=& \beta_2 T_2^m(x) + \alpha - \i^2_m(x) -(\beta_1 T_1^m(x) +
\alpha - \i^1_m(x)) \\ \nonumber &=& \beta_2(T_2^m(x) - T_1^m(x)) + \delta T_1^m(x) \geq
\beta_2^m \delta x. \eea On the other hand, $1 \geq T_2^m(x) - T_1^m(x) \geq \beta_2^{m-1}
\delta x$. Thus $\delta \leq \frac{\beta_2^{-m+1}}{x}$ for all $m$ such that $\i^1_{[0,m)}
= \i^2_{[0,m)}$. If $x=0$, then $T_1(x) = T_2(x) = \alpha$ and we can apply the first
statement to $y = \alpha > 0$. $\Box$

\begin{prop} \label{prop uv to h(Suv) continuous}
Let $(\uu,\uv),(\uu',\uv') \in \Sk^2$ satisfy \eqref{eq uu uv}.
Let $L \in \N$ and suppose that $\uu$, $\uu'$ have a common prefix
of length larger than $L$ and $\uv$, $\uv'$ have a common prefix
of length larger than $L$. Then for all $\delta>0$ there exists
$L(\delta)$ such that for any $L\geq L(\delta)$,
\be \nonumber
|\htop(\Sigma_{\uu',\uv'})-\htop(\Suv)| \leq \delta\,.
\ee
\end{prop}

This proposition is a stronger reformulation of  Proposition
9.3.15 in \cite{BrBr04}. It follows from the proof of Proposition
9.3.15 given in this book, except that the argument at the very
end of the proof is incomplete; but it is completed in
\cite{BrTo08}.  Now we can state our first theorem
and his corollary about the normality of orbits under $\Tab$. The
proof of the theorem is inspired by the proof of Theorem C in
\cite{Sc97}, where the case $x=1$ and $\alpha=0$ is considered.

\begin{theo} \label{theo selfnormality for iab alpha fixed}
Let $x \in [0,1)$ and $\alpha \in [0,1)$ excepted $(x,\alpha) = (0,0)$. Then the set \be
\nonumber \{\beta>1: \text{the orbit of $\iab(x)$ under $\sigma$ is $\hmuab$-normal}\} \ee
has full $\lambda$-measure.
\end{theo}

\begin{cl} \label{cl selfnormality for Tab alpha fixed}
Let $x \in [0,1)$ and $\alpha \in [0,1)$ excepted $(x,\alpha) = (0,0)$. Then the set \be
\nonumber \{\beta>1: \text{the orbit of $x$ under $\Tab$ is $\muab$-normal}\} \ee has full
$\lambda$-measure.
\end{cl}

Remark that the theorem and its corollary may also be formulated for $x \in (0,1]$ using a
left-continuous extension of $\Tab$ on $(0,1]$ and a coding $\iab$ defined using intervals
$(a_j,a_{j+1}]$ for all $j \in \A$.

\textbf{Proof of the theorem:} We briefly sketch the proof. We use
the uniqueness of the measure of maximal entropy $\hmuab$: for
$\ux \in \Sab$ not $\hmuab$-normal, there exists $\nu \in
V_\sigma(\ux)$ such that $h(\nu) < h(\hmuab) = \log\beta$. The
main idea is to imbed $\{\iab(x): \beta \in [\beta_1,\beta_2]\}$
in a shift space $\Sigma^* := \Sigma_{\uu^*,\uv^*}$ with $\uu^*$
and $\uv^*$ well chosen. Writing $D^* \subset \Sigma^*$ for the
range of the imbedding, we estimate the Hausdorff dimension of the
subset of $D^*$ corresponding to points $\iab(x)$ which are not
$\hmuab$-normal. Then we estimate the coefficient of H\"older
continuity of the application $\rho_*$ defined as the inverse of
the imbedding. This gives us an estimate of the Hausdorff
dimension of the non $\hmuab$-normal points in the interval
$[\beta_1,\beta_2]$.

To obtain uniform estimates, we restrict our proof to the interval $[\ub,\ob]$ with $1 <
\ub < \ob < \infty$. This is sufficient, since there exist a countable cover of
$(1,\infty)$ with such intervals. Let $k := \lceil \alpha + \ob \rceil$ and $\Omega :=
\{\beta\in [\ub,\ob]: \iab(x) \text{ is not $\hmuab$-normal}\}$. For $\beta \in \Omega$, we
have $V_\sigma(\iab(x)) \neq \{\hmuab\}$. Since $\hmuab$ is the unique measure of maximal
entropy $\log\beta$, there exist $N \in \N$ and $\nu \in V_\sigma(\iab(x))$ such that
$h(\nu) < (1-1/N)\log\beta$. Setting \be \nonumber \Omega_N := \{\beta \in [\ub,\ob]:
\exists \nu \in V_\sigma(\iab(x)) \text{ s.t. } h(\nu) < (1-1/N) \log\beta\}, \ee we have
$\Omega = \bigcup_{N\geq1} \Omega_N$. We will prove that $\dim_H \Omega_N <1$, so that
$\lambda(\Omega_N) = 0$ for all $N\geq1$.

For $N \in \N$ fixed, define $\eps:= \frac{\ub \log\ub}{2N-1}>0$ and
$\delta:=\log\big(1+\eps/\ob\big)$. Choose  $L\geq L(\delta)$ (Proposition \ref{prop uv to
h(Suv) continuous}). Consider the family of subsets of $[\ub,\ob]$ of the following type
\be \nonumber J(\uw,\uw')  = \{\beta \in [\ub,\ob]: \uab_{[0,L)} = \uw, \vab_{[0,L)} =
\uw'\} \ee where $\uw,\uw'$ are two words of length $L$. $J(\uw,\uw')$ is either empty or
it is an interval, since the applications $\beta \mapsto \uab$ and $\beta \mapsto \vab$ are
both monotone increasing. Moreover, $[\ub,\ob] = \bigcup_{\uw,\uw'} J(\uw,\uw')$ where the
union is finite, since the set of words of length $L$ in $\As$ has finite cardinality. We
want to work with closed intervals, thus we cover the non-closed $J(\uw,\uw')$ with
countably many closed intervals if necessary. For example, if $J(\uw,\uw') = (a,b]$, we
write $J(\uw,\uw') = \bigcup_{m\geq1}[a+1/m,b]$. We prove that $\lambda(\Omega_N \cap
[\beta_1,\beta_2]) = 0$ where $\beta_1<\beta_2$ are such that $\uu^{\ab_1}_{[0,L)} =
\uu^{\ab_2}_{[0,L)}$ and $\uv^{\ab_1}_{[0,L)} = \uv^{\ab_2}_{[0,L)}$.

Let $\uu^j = \uu^{\ab_j}$ and $\uv^j = \uv^{\ab_j}$. Using \eqref{eq uab vab} and the
monotonicity of $\beta \mapsto \uab$ and $\beta \mapsto \vab$, we have \be
\nonumber \begin{array}{l} \uu^1 \preceq \sigma^n \uu^1 \preceq \uv^1 \preceq \uv^2 \\
\uu^1 \preceq \uu^2 \preceq \sigma^n \uv^2 \preceq \uv^2 \end{array} \quad \forall n \geq0.
\ee Hence the couple $(\uu^1,\uv^2)$ satisfy \eqref{eq uu uv} and we set $\Sigma^* :=
\Sigma_{\uu^1,\uv^2}$ and \be \nonumber D^* := \{\uz \in \Sigma^*:
\exists \beta \in [\beta_1,\beta_2] \text{ s.t. } \uz = \iab(x)\}.
\ee We define an application $\rho_*: D^* \to [\beta_1,\beta_2]$
by $\rho_*(\uz) = \beta \Leftrightarrow \iab(x) = \uz$. This
application is well defined: by definition of $D^*$, for all $\uz
\in D^*$ there exists a $\beta$ such that $\uz = \iab(x)$;
moreover this $\beta$ is unique, since by Lemma \ref{lm separation
of coding in Tab}, $\beta \mapsto \iab(x)$ is strictly increasing.
On the other hand, for all $\beta \in [\beta_1,\beta_2]$, we have
from \eqref{eq Sab} \be \nonumber \uu^1 \leq \uab \leq \sigma^n
\iab(x) \leq \vab \leq \uv^2 \quad \forall n\geq0, \ee whence
$\iab(x) \in \Sigma^*$ and $\rho_*: D^* \to [\beta_1,\beta_2]$ is
surjective. Let $\beta_* := \htop(\Sigma^*)$;  then by Proposition
\ref{prop uv to h(Suv) continuous}
\be \label{eq
bstar near b1}
\beta_* -\beta_1\leq{\rm e}^{\log\beta_1}\big({\rm
e}^{\delta}-1\big)\leq\eps\,.
\ee
Let us compute the coefficient of H\"older continuity of $\rho_*:
(D^*,d_{\beta_*}) \to [\beta_1,\beta_2]$. Let $\uz \neq \uz' \in
D^*$ and $n = \min\{l\geq0: z_l \neq z_l'\}$, then
$d_{\beta_*}(\uz,\uz') = \beta_*^{-n}$. By Lemma \ref{lm
separation of coding in Tab}, there exists $C$ such that \be
\nonumber |\rho_*(\uz)-\rho_*(\uz')| \leq C \rho_*(\uz)^{-n} \leq
C \beta_1^{-n} =
C(d_{\beta_*}(\uz,\uz'))^{\frac{\log\beta_1}{\log\beta_*}}. \ee We
may choose $C$ independently of $\beta$, since we work on the
compact interval $[\beta_1,\beta_2] \subset (1,\infty)$. By
equation \eqref{eq bstar near b1} and the choice of $\eps$, we
have \bea \nonumber \beta_* - \beta_1 \leq \frac{\ub \log
\ub}{2N-1} &\Rightarrow& \beta_* - \beta_1 \leq \frac{\beta_1 \log
\beta_1}{2N-1} \\ \nonumber &\Leftrightarrow& 1 + \frac{\beta_* -
\beta_1}{\beta_1 \log \beta_1} \leq 1 + \frac{1}{2N-1} \\ \nonumber &\Leftrightarrow&
\frac{\log\beta_1 + \frac{\beta_* - \beta_1}{\beta_1}} {\log \beta_1}  \leq \frac{2N}{2N-1}
\\ \nonumber &\Rightarrow& \frac{\log \beta_1}{\log \beta_*} \geq \frac{\log \beta_1}{\log
\beta_1 + \frac{\beta_* - \beta_1}{\beta_1}} \geq 1 - \frac{1}{2N}. \eea In last line, we
use the concavity of the logarithm, so the first order Taylor
development is an upper estimate. Thus $\rho_*$ has
H\"older-exponent  $1-\frac{1}{2N}$.

Define \be \nonumber G_N^* := \{\uz \in \Sigma^*: \exists \nu \in V_\sigma(\uz) \text{ s.t.
} h(\nu) < (1-1/N) \log\beta_*\}. \ee Let $\beta \in \Omega_N \cap [\beta_1,\beta_2]$. Then
there exists $\nu \in V_\sigma(\iab(x))$ such that \be \nonumber h(\nu) < (1-1/N)\log\beta
\leq (1-1/N) \log\beta_*. \ee Since $\iab(x) \in D^* \subset \Sigma^*$, we have $\iab(x)
\in G_N^*$. Using the surjectivity of $\rho_*$, we obtain $\Omega_N \cap [\beta_1,\beta_2]
\subset \rho_*(G_N^* \cap D^*)$. We claim that $\htop(G_N^*) \leq (1-1/N)\log\beta_*$. This
implies, using Lemmas \ref{lm dimH Holder continuity} and \ref{lm htop dimH}, \be \nonumber
\dim_H(\Omega_N \cap [\beta_1,\beta_2]) \leq \dim_H \rho_*(G_N^* \cap D^*) \leq
\frac{\dim_H G_N^*}{1-\frac{1}{2N}} \leq \frac{\htop(G_N^*)} {(1-\frac{1}{2N})\log\beta_*}
\leq \frac{1-\frac{1}{N}} {1-\frac{1}{2N}} < 1. \ee Thus $\lambda(\Omega_N \cap
[\beta_1,\beta_2]) = 0$.

It remains to prove $\htop(G_N^*) \leq (1-1/N)\log\beta_*$. Recall that $h(\nu) = \lim_n
\frac{1}{n} H_n(\nu)$, where $H_n(\nu)$ is the entropy of $\nu$ with respect to the algebra
of cylinder sets of length $n$ and that $\frac{1}{n} H_n(\nu)$ is decreasing. For all
$m\geq1$, we set \bea \nonumber F_N^*(m) &:=& \{\nu \in M(\Sigma^*,\sigma): \frac{1}{m}
H_m(\nu) \leq (1-1/N)\log\beta_*\} \\ \nonumber G_N^*(m) &:=& \{\uz \in \Sigma^*:
V_\sigma(\uz) \cap F_N^*(m) \neq \emptyset\}. \eea Let $\uz \in G_N^*$, then there exists
$\nu \in V_\sigma(\uz)$ such that $h(\nu) < (1-\frac{1}{N})\log\beta_*$. Since
$\frac{1}{m}H_m(\nu) \downarrow h(\nu)$, there exists $m\geq1$ such that $\frac{1}{m}
H_m(\nu) \leq (1-1/N)\log\beta_*$, whence $\nu \in F_n^*(m)$ and $\uz \in G_N^*(m)$. This
implies $G_N^* \subset \bigcup_{m\geq1} G_N^*(m)$. Since $H_m(\cdot)$ is continuous,
$F_N^*(m)$ is closed for all $m\geq1$. Finally we obtain using Theorem \ref{theo htop ^K G}
\be \nonumber \htop(G_N^*) = \sup_m \htop(G_N^*(m)) \leq \sup_m \sup_{\nu \in F_N^*(m)}
h(\nu) \leq \sup_m \sup_{\nu \in F_N^*(m)} \frac{1}{m} H_m(\nu) \leq (1-1/N)\log\beta_*.
\quad \Box \ee

\textbf{Proof of the Corollary:} Let $\beta>1$ be such that the orbit of $\iab(x)$ under
$\sigma$ is $\hmuab$-normal. Let $f \in C([0,1])$, then $\hat{f}: \Sab \to \R$ defined by
$\hat{f} := f \circ \phiab$ is continuous, since $\phiab$ is continuous. Using $\muab :=
\hmuab \circ (\phiab)^{-1}$, we have \bea \nonumber \int_{[0,1]} f d\muab &=& \int_{\Sab}
\hat{f} d\hmuab = \lim_{n\to\infty} \sum_{i=0}^{n-1} \hat{f}(\sigma^i \iab(x)) \\ \nonumber
&=& \lim_{n\to\infty} \sum_{i=0}^{n-1} f(\phiab(\sigma^i \iab(x))) = \lim_{n\to\infty}
\sum_{i=0}^{n-1} f(\Tab^i(x)). \eea The second equality comes from the $\hmuab$-normality
of the orbit of $\iab(x)$ under $\sigma$, the last one is \eqref{eq T^n = phi sigma^n i}
which is true for all $x \in [0,1)$ with our convention for the extension of $\Tab$ and
$\iab$ on $[0,1)$. $\Box$

The next step is to consider the question of $\muab$-normality in the whole plane $(\ab)$
instead of working with $\alpha$ fixed. Define $\cR := [0,1) \times (1,\infty)$.

\begin{theo} \label{theo selfnormality for Tab}
For all $x \in [0,1)$, the set \be \nonumber \cN(x) := \{(\alpha,\beta) \in \cR: \text{the
orbit of $x$ under $\Tab$ is $\muab$-normal}\} \ee has full $\lambda^2$-measure.
\end{theo}

\Pr We have only to prove that $\cN(x)$ is measurable and to apply Fubini's Theorem and
Corollary \ref{cl selfnormality for Tab alpha fixed}. The first step is to prove that for
all $x\in[0,1)$ and all $n\geq0$, the applications $(\ab) \mapsto \iab(x)$ and $(\ab)
\mapsto \Tab^n(x)$ are measurable. First remark that for all $n\geq1$ \be \label{eq
Tab^n(x)} \Tab^n(x) = \beta^n x + \alpha \frac{\beta^n-1}{\beta-1} - \sum_{j=0}^{n-1}
\iab_j(x) \ \beta^{n-j-1}. \ee The proof by induction is immediate. To prove that $(\ab)
\mapsto \iab(x)$ is measurable, it is enough to prove that for all $n\geq0$ and for all
words $\uw \in \As$ of length $n$ \be \nonumber \{(\ab) \in \cR: \iab_{[0,n)}(x) = \uw \}
\ee is measurable, since the $\sigma$-algebra on $\Sk$ is generated by the cylinders. This
set is the subset of $\R^2$ such that \be \nonumber
\begin{cases} \beta>1
\\ 0 \leq \alpha < 1 \\ w_j < \beta \Tab^j(x) + \alpha \leq w_j+1 & \forall 0\leq j < n
\end{cases} \ee Using \eqref{eq Tab^n(x)}, this system of
inequations can be rewritten \be \nonumber \begin{cases} \beta>1 \\ 0 \leq  \alpha < 1 \\
\alpha > \frac{\beta-1} {\beta^{j+1}-1} \left(\sum_{i=0}^j w_i \beta^{j-i} - \beta^{j+1}x
\right) & \forall 0\leq j < n \\ \alpha \leq \frac{\beta-1} {\beta^{j+1}-1} \left(1 +
\sum_{i=0}^j w_i \beta^{j-i} - \beta^{j+1}x \right) & \forall 0\leq j < n \end{cases} \ee
From this, the measurability of $\iab$ follows. If $(\ab) \mapsto \iab(x)$ is measurable,
then by formula \eqref{eq Tab^n(x)}, $(\ab) \mapsto \Tab^n(x)$ is clearly measurable for
all $n\geq0$. Then for all $f \in C([0,1])$ and all $n\geq1$, the application $(\ab)
\mapsto S_n(f) := \frac{1}{n} \sum_{i=0}^{n-1} f(\Tab^i(x))$ is measurable and consequently
\be \nonumber \{(\ab): \lim_{n\to\infty} S_n(f) \text{ exists} \} \ee is a measurable set.

On the other hand, if $f \in C([0,1])$, then $(\ab) \mapsto \int f d\muab$ is measurable.
Indeed \be \nonumber \int f d\muab = \int f h_{\ab} d\lambda \ee and in view of equation
\eqref{eq density muab} and the measurability of $(\ab) \mapsto \Tab(x)$, the application
$(\ab) \mapsto h_{\ab}$ is clearly measurable. Therefore \be \nonumber \{(\ab):
\lim_{n\to\infty} S_n(f) = \int f d\muab \} \ee is measurable for all $f \in C([0,1])$. Let
$\{f_m\}_{m\in\N} \subset C([0,1])$ be countable subset which is dense with respect to the
uniform convergence. Then setting \be \nonumber D_m := \{(\ab) \in \cR: \lim_{n\to\infty}
S_n(f_m) = \int f_m d\muab \}, \ee we have $\cN(x) = \bigcap_{m\in\N} D_m$, whence it is a
measurable set. $\Box$

We have shown that for a given $x \in [0,1)$, the orbit of $x$ under $\Tab$ is
$\muab$-normal for almost all $(\ab)$. The orbits of $0$ and $1$ are of particular interest
(see equations \eqref{eq density muab} or \eqref{eq Sab}). Now we show that by any point
$(\alpha_0,\beta_0)$, there passes a curve defined by $\alpha = \alpha(\beta)$ such that
the orbit of $0$ under $\Tabb$ is $\muabb$-normal for at most one $\beta$. A trivial
example of such a curve is $\alpha =0$, since $x=0$ is a fixed point. The idea is to
consider curves along which the coding of $0$ is constant, ie to define $\alpha(\beta)$
such that $\uabb$ is constant. The results below depend on reference \cite{FaPf08}, where
we solve the following inverse problem: given $\uu$ and $\uv$ verifying \eqref{eq uu uv},
can we find $\ab$ such that $\uu = \uab$ and $\uv = \vab$ ?

Let \be \nonumber \cU := \{\uu: \exists \ (\ab) \in \cR \text{
s.t. } \uu = \uab \}. \ee We define an equivalence relation in
$\cR$ by \be \nonumber (\ab) \sim (\alpha',\beta') \iff \uab =
\uu^{\alpha',\beta'}. \ee An equivalence class is denoted by
$[\uu]$. The next lemma describes $[\uu]$.

\begin{lm} \label{lm [uu] is a continuous curve}
Let $\uu \in \cU$ and set \be \nonumber \alpha(\beta) = (\beta-1) \sum_{j\geq0}
\frac{u_j}{\beta^{j+1}}. \ee Then there exists $\beta_{\uu} \geq1$ such that \be \nonumber
[\uu] = \{(\abb): \beta \in I_{\uu} \} \ee with $I_{\uu} = (\beta_{\uu}, \infty)$ or
$I_{\uu} = [\beta_{\uu},\infty)$.
\end{lm}

\Pr If $\uu = 000\dots$, then the statement is trivially true with
$\alpha(\beta) \equiv 0$ and $\beta_{\uu} = 1$. Suppose $\uu \neq
000\dots$. First we prove that \be \nonumber (\ab) \sim
(\alpha',\beta) \implies \alpha = \alpha' \ee then \be \nonumber
(\ab) \in [\uu] \implies (\alpha(\beta'),\beta') \in [\uu] \quad
\forall \beta'\geq\beta. \ee Let $(\ab) \in [\uu]$. Using
\eqref{eq T^n = phi sigma^n i}, we have $\phiab(\sigma \uu) =
\Tab(0) = \alpha$. Since the map $\alpha \mapsto \phiab(\sigma
\uu) - \alpha$ is continuous and strictly decreasing (Lemmas 3.5
and 3.6 in \cite{FaPf08}), the first statement is true. Let
$\beta'>\beta$. By Corollary 3.1 in \cite{FaPf08}, we have that
$\phiab(\sigma \uu) > \phi^{\ab'}(\sigma \uu)$. Therefore there
exists a unique $\alpha'<\alpha$ such that
$\phi^{\alpha',\beta'}(\sigma \uu) = \alpha'$. We prove that
$\uu^{\alpha',\beta'} = \uu$. By point 1 of Proposition 2.5 in
\cite{FaPf08}, we have $\uu \preceq \uu^{\alpha',\beta'}$. By
Proposition 3.3 in \cite{FaPf08}, we have \be \nonumber
\htop(\Sigma_{\uu,\uv^{\alpha',\beta'}}) =
\htop(\Sigma_{\alpha',\beta'}) = \log\beta'. \ee Since $\Sab =
\Sigma_{\uu,\vab}$ and $\beta'>\beta$, we must have $\vab \prec
\uv^{\alpha',\beta'}$. Therefore
\be \nonumber
\begin{cases}
\uu
\preceq \sigma^n \uu \prec \vab \prec \uv^{\alpha',\beta'} \\
\uu\preceq \uu^{\alpha',\beta'} \prec \sigma^n
\uv^{\alpha',\beta'} \preceq \uv^{\alpha',\beta'}
\end{cases}
\quad \forall n\geq0,
\ee
are the inequalities (4.1) in \cite{FaPf08} for the pair
$(\uu,\uv^{\alpha',\beta'})$. We can apply Proposition 3.2 and
Theorem 4.1 in \cite{FaPf08} to this pair and get $\uu =
\uu^{\alpha',\beta'}$. It remains to show that $\alpha' =
\alpha(\beta')$. Following the definition of the
$\varphi$-expansion of R\'enyi, we have for all $x \in [0,1)$ and
all $n\geq0$
\be \nonumber x = \sum_{j=0}^{n-1} \frac{\iab_j(x)-\alpha}{\beta^{j+1}} +
\frac{\Tab^n(x)}{\beta^n}. \ee Since $\Tab^n(x) \in [0,1)$, for all $\beta>1$ we find an
explicit expression for $\phiab$ on $\Sab$ \be \nonumber x = \sum_{j\geq0}
\frac{\iab_j(x)-\alpha}{\beta^{j+1}}. \ee In particular, applying this equation to $x=0$,
we have for all $(\ab) \in \cR$ \be \nonumber \alpha = (\beta-1) \sum_{j\geq0}
\frac{u^{\ab}_j}{\beta^{j+1}}. \ee Since for all $\beta>\beta_{\uu}$, we have $\uu \in
\Sab$, this complete the proof. $\Box$

For each $\uu \in \cU$, the equivalence class $[\uu]$ defines an
analytic curve in $\cR$, which is strictly monotone decreasing
(excepted for $\uu = 000\dots$),
\be \nonumber
[\uu] = \{(\ab): \alpha
= (\beta-1) \sum_{j\geq0} \frac{u_j}{\beta^{j+1}}, \ \beta \in
I_{\uu} \}.
\ee
There curves are disjoint two by two and their
union is $\cR$.

\begin{theo} \label{theo selfnormality for Tab along curves}
Let $(\ab) \in \cR$, $\uu = \uab$ and define $\alpha(\beta)$ and $\beta_{\uu}$ as in Lemma
\ref{lm [uu] is a continuous curve}. Then for all $\beta>\beta_{\uu}$, the orbit of $x=0$
under $\Tabb$ is not $\muabb$-normal.
\end{theo}

\Pr Let $\hat{\nu} \in M(\Sigma_k,\sigma)$ (with $k$ large enough)
be a cluster point of $\{\cE_n(\uu)\}_{n\geq1}$. By Lemma \ref{lm
[uu] is a continuous curve}, $\uabb = \uu$ for any
$\beta>\beta_{\uu}$. Therefore
\be\nonumber
h(\hat{\nu}) \leq \htop(\Sigma_{\abb}) = \log\beta \qquad \forall
\beta>\beta_{\uu}
\ee
and $\hat{\nu}$ is not a measure of maximal entropy, as well as
$\nu_{\beta} := \hat{\nu} \circ (\phiabb)^{-1}$ \cite{Ho79}, for
all $\beta>\beta_{\uu}$.
$\Box$

Recall that \be \nonumber \cN(0) = \{(\ab) \in \cR: \text{the orbit of $0$ under $\Tab$ is
$\muab$-normal} \}. \ee By Theorem \ref{theo selfnormality for Tab}, $\cN(0)$ has full
Lebesgue measure. On the other hand, by Theorem \ref{theo selfnormality for Tab along
curves}, we can decompose $\cR$ into a family of disjoint analytic curves such that each
curve meets $\cN(0)$ in at most one point. This situation is very similar to the one
presented in \cite{Mi97} by Milnor following an idea of Katok.

\section{Normality in generalized $\beta$-maps} \label{sec generalized beta-maps}

In this section, we consider another class of piecewise monotone continuous applications,
the generalized $\beta$-maps. Introduced by G\'ora in \cite{Go07}, they have only one
critical orbit like $\beta$-maps, but they admit increasing and decreasing laps. A family
$\{\Tb\}_{\beta>1}$ of generalized $\beta$-maps is defined by $k\geq2$ and a sequence
$s=(s_n)_{0\leq n<k}$ with $s_i \in \{-1,1\}$. For any $\beta \in (k-1,k]$, let $a_j =
j/\beta$ for $j = 0,\dots,k-1$ and $a_k = 1$. Then for all $j= 0,\dots,k-1$, the map $f_j=
I_j \to [0,1]$ is defined by \be \nonumber f_j(x) :=
\begin{cases} \beta x \mod 1 & \text{if } s_j = +1 \\ 1 -(\beta x \mod 1) & \text{if } s_j
= -1. \end{cases} \ee In particular when $s = (1,-1)$, then $\Tb$ is a tent map. Here we
left the map undefined on $a_j$ for $j = 1,\dots,k-1$.

G\'ora constructed the unique measure $\mub$ absolutely continuous with respect to Lebesgue
measure (Theorem 6 and Proposition 8 in \cite{Go07}). Using the same argument as Hofbauer
in \cite{Ho79bis}, we deduce that a measure of maximal entropy is always absolutely
continuous with respect to Lebesgue measure, hence the measure $\mub$ is the unique measure
of maximal entropy. Let $k = \lceil \beta \rceil$ and let us denote $\ib$ for the coding
map under $\Tb$, $\phib := (\ib)^{-1}$ for the inverse of the coding map, $\Sb :=
\Sigma_{\Tb}$ and $\eb := \lim_{x\uparrow1} \ib(x)$. Now it is easy to check that formula
\eqref{eq Sigma_T} becomes \be \label{eq Sb} \Sb = \{\ux \in \Sk: \sigma^n \ux \preceq \eb
\quad \forall n\geq0\} \ee and inequations \eqref{eq uuj uvj} become \be \label{eq ep}
\sigma^n \eb \preceq \eb \quad \forall n\geq0. \ee It is known that the dynamical system
$(\Sb,\sigma)$ has topological entropy $\log\beta$ and, by general works of Hofbauer in
\cite{Ho79}, it has a unique measure of maximal entropy $\hmub$ such that $\mub = \hmub
\circ (\phib)^{-1}$.

As in the previous section, we state two lemmas which we need for the proof of the main
theorem of this section. We study the normality only of $x=1$, so these lemmas are
formulated only for $x=1$. Let $S_n(\beta) \equiv S_n$ and $S(\beta) \equiv S$ be defined
by \eqref{eq X_n S_n S}.

\begin{lm} \label{lm 1 notin S for almost all beta}
For any family of generalized $\beta$-maps defined by $(s_n)_{0\leq n<k}$, the set
$\{\beta\in (k-1,k]: 1 \in S(\beta) \}$ is countable.
\end{lm}

\Pr For a fixed $n\geq1$, we study the map $\beta \mapsto \Tb^n(1)$. This map is well
defined everywhere in $(k-1,k]$ excepted for finitely many points and it is continuous on
each interval where it is well defined. Indeed this is true for $n=1$. Suppose it is true
for $n$, then $\Tb^{n+1}(1)$ is well defined and continuous wherever $\Tb^n(1)$ is well
defined and continuous, excepted when $\Tb^n(1) \in S_0(\beta)$. By the induction
hypothesis, there exists a finite family of disjoint open intervals $J_i$ and continuous
functions $g_i: J_i \to [0,1]$ such that $(k-1,k] \bk (\bigcup_i J_i)$ is finite and \be
\nonumber \Tb^n(x) = g_i(\beta) \quad \text{if } \beta \in J_i. \ee Then \be \nonumber
\{\beta \in (k-1,k]: \Tb^n(1) \text{ is well defined and } \Tb^n(1) \in S_0(\beta) \} =
\bigcup_{i,j} \{\beta \in J_i: g_i(\beta) = \frac{j}{\beta} \}. \ee We claim that $\{\beta
\in J_i: g_i(\beta) = \frac{j}{\beta} \}$ has finitely many points. From the form of the
map $\Tb$, it follows immediately that each $g_i(\beta)$ is a polynomial of degree $n$.
Since $\beta>1$, \be \nonumber g_i(\beta) = \frac{j}{\beta} \quad \iff \quad \beta
g_i(\beta) -j = 0. \ee This polynomial equation has at most $n+1$ roots. In fact, using the
monotonicity of the map $\beta \mapsto \eb$, we can prove that this set has at most one
point. The lemma follows, since $S(\beta) = \bigcup_{n\geq0} S_n(\beta)$. $\Box$

\begin{lm} \label{lm separation of coding in Tb}
Consider a family $\{\Tb\}_{\beta>1}$ of generalized $\beta$-maps defined by a sequence $s
= (s_n)_{n\geq0}$. Let $1<\beta_1\leq\beta_2$ and define $l = \min\{n\geq0: \ue^1_n \neq
\ue^2_n\}$ with $\ue^j = \ue^{\beta_j}$ for $j=1,2$. \newline If $k\geq3$, for all
$\beta_0>2$, there exists $K$ such that $\beta_1 \geq \beta_0$ implies \be \nonumber
\beta_2 - \beta_1 \leq K \beta_2^{-l}. \ee \newline If $s = (+1,+1)$, then \be \nonumber
\beta_2-\beta_1 \leq \beta_2^{-l+1}. \ee If $s = (+1,-1)$ or $(-1,+1)$, then for all
$\beta_0>1$, there exists $K$ such that $\beta_1 \geq \beta_0$ implies \be \nonumber
\beta_2-\beta_1 \leq K \beta_2^{-l}. \ee If $s = (-1,-1)$, then there exists $\beta_0>1$
and $K$ such that $\beta_1 \geq \beta_0$ implies \be \nonumber \beta_2-\beta_1 \leq K
\beta_2^{-l}. \ee
\end{lm}

The proof is very similar to the proof of Brucks and Misiurewicz for Proposition 1 of
\cite{BrMi96}, see also Lemma 23 of Sands in \cite{Sa93}.

\Pr Let $\delta := \beta_2-\beta_1\geq0$ and denote $T_j = T_{\beta_j}$ and $\i^j =
\i^{\beta_j}$ for $j=1,2$. Let $a_1, a_2 \in [0,1]$ such that $r := \i^1_0(a_1) =
\i^2_0(a_2)$. Considering four cases according to the signs of $a_2-a_1$ and $s_r$, we have
\be \nonumber |T_2(a_2)-T_1(a_1)| \geq \beta_2|a_2-a_1| - \delta. \ee Applying $n$ times
this formula, we find that $\i^1_{[0,n)}(a_1) = \i^2_{[0,n)}(a_2)$ implies \be \nonumber
|T_2^n(a_2) - T_1^n(a_1)| \geq \beta_2^n \left( |a_2-a_1| - \frac{\delta}{\beta_2-1}
\right). \ee Consider the case $k\geq3$. Then $a_i = T_i(1)$ for $i = 1,2$ are such that
$|a_2-a_1| = \delta > \frac{\delta}{\beta_0-1} \geq \frac{\delta}{\beta_2-1}$. Using
$|T_2^n(a_2) - T_1^n(a_1)| \leq 1$, we conclude that for all $\beta_0 \leq \beta_1 \leq
\beta_2$, if $\ue^1_{[0,n)} = \ue^2_{[0,n)}$ then \be \nonumber \delta \leq
\frac{\beta_0-1}{\beta_0-2} \ \beta_2^{-n+1}. \ee For the case $s = (+1,+1)$, we can apply
Lemma \ref{lm separation of coding in Tab} with $\alpha=0$ and $x=1$. \newline The case $s
= (+1,-1)$ or $(-1,+1)$ is considered in Lemma 23 of \cite{Sa93}. \newline For the case $s
= (-1,-1)$: for a fixed $n$, we want to find $\beta_0$ such that for all $\beta_0 \leq
\beta_1 \leq \beta_2$ we have \be \label{eq lm separation in Tb} |T_2^n(1) - T_1^n(1)| >
\frac{\delta}{\beta_2-1}. \ee Then we conclude as in the case $k\geq3$. The formula
\eqref{eq lm separation in Tb} is true, if $|\frac{d}{d\beta} \Tb^n(1)| >
\frac{1}{\beta-1}$ for all $\beta \geq \beta_0$. When $n$ increases, $\beta_0$ decreases.
With $n = 3$, we have $\beta_0 \approx 1.53.$ $\Box$

In the tent map case, the separation of orbits is proved for $\beta \in (\sqrt{2},2]$ and
then extended arbitrarily near $\beta_0=1$ using the renormalization. In the case $s =
(-1,-1)$, there is no such argument and we are forced to increase $n$ to obtain a lower
bound $\beta_0$. With the help of a computer, we obtain $\beta_0 \approx 1.27$ for $n=12$.
For more details, see \cite{Fa08}.

Now we turn to the question of normality for generalized $\beta$-maps. The structure of the
proof is very similar to the proof of Theorem \ref{theo selfnormality for iab alpha fixed}
and Corollary \ref{cl selfnormality for Tab alpha fixed}.

\begin{theo} \label{theo selfnormality for Tb}
Consider a family $\{\Tb\}_{k-1<\beta\leq k}$ of generalized $\beta$-maps defined by a
sequence $s = (s_n)_{0\leq n<k}$. Let $\beta_0$ be defined as in Lemma \ref{lm separation
of coding in Tb} according to $s$. Then the set \be \nonumber \{\beta>\beta_0: \text{the
orbit of $\eb$ under $\sigma$ is $\hmub$-normal}\} \ee has full $\lambda$-measure.
\end{theo}

\begin{cl} \label{cl selfnormality for Tb}
Consider a family $\{\Tb\}_{\beta>1}$ of generalized $\beta$-maps defined by a sequence $s
= (s_n)_{n\geq0}$. Let $\beta_0$ be defined as in Lemma \ref{lm separation of coding in Tb}
according to $s$. Then the set \be \nonumber \{\beta>\beta_0: \text{the orbit of $1$ under
$\Tb$ is $\mub$-normal}\} \ee has full $\lambda$-measure.
\end{cl}

\textbf{Proof of Theorem:} Let \be \nonumber B_0 := \{\beta \in (\beta_0,\infty): 1 \notin
S(\beta)\}. \ee From Lemma \ref{lm 1 notin S for almost all beta}, this subset has full
Lebesgue measure. To obtain uniform estimates, we restrict our proof to the interval
$[\ub,\ob]$ with $\beta_0 < \ub < \ob < \infty$. Let $k := \lceil \ob \rceil$ and $\Omega
:= \{\beta\in [\ub,\ob]\cap B_0: \eb \text{ is not $\hmub$-normal}\}$. As before, setting
\be \nonumber \Omega_N := \{\beta \in [\ub,\ob] \cap B_0: \exists \nu \in V_\sigma(\eb)
\text{ s.t. } h(\nu) < (1-1/N) \log\beta\}, \ee we have $\Omega = \bigcup_{N\geq1}
\Omega_N$. We prove that $\dim_H \Omega_N <1$. For $N \in \N$ fixed, define $\eps:=
\frac{\ub \log\ub}{2N-1}>0$ and $L$ such that $\eb_{[0,L)} = \ue^{\beta'}_{[0,L)}$ implies
$|\beta-\beta'| \leq \eps$ (see Lemma \ref{lm separation of coding in Tb}). Consider the
family of subsets of $[\ub,\ob]$ of the following type \be \nonumber J(\uw)  = \{\beta \in
[\ub,\ob]: \eb_{[0,L)} = \uw\} \ee where $\uw$ is a word of length $L$. $J(\uw)$ is either
empty or it is an interval. We cover the non-closed $J(\uw)$ with countably many closed
intervals if necessary. We prove that $\lambda(\Omega_N \cap [\beta_1,\beta_2]) = 0$ where
$\beta_1<\beta_2$ are such that $\ue^{\beta_1}_{[0,L)} = \ue^{\beta_2}_{[0,L)}$.

Let $\ue^j = \ue^{\beta_j}$. Let \be \nonumber D^* := \{\uz \in \Sigma_{\ue^2}: \exists
\beta \in [\beta_1,\beta_2]\cap B_0 \text{ s.t. } \uz = \eb \}. \ee Define $\rho_*: D^* \to
[\beta_1,\beta_2]\cap B_0$ by $\rho_*(\uz) = \beta \Leftrightarrow
\eb = \uz$. As before, from formula \eqref{eq Sb} and strict
monotonicity of $\beta \mapsto \eb$, we deduce that $\rho_*$ is
well defined and surjective. We compute the coefficient of
H\"older continuity of $\rho_*: (D^*,d_{\beta_*}) \to
[\beta_1,\beta_2]$. Let $\uz \neq \uz' \in D^*$ and $n =
\min\{l\geq0: z_l \neq z_l'\}$, then $d_{\beta_*}(\uz,\uz') =
\beta_*^{-n}$. By Lemma \ref{lm separation of coding in Tb}, there
exists $C$ such that
\be \nonumber |\rho_*(\uz)-\rho_*(\uz')| \leq C \rho_*(\uz)^{-n}
\leq C \beta_1^{-n} =
C(d_{\beta_*}(\uz,\uz'))^{\frac{\log\beta_1}{\log\beta_*}}. \ee By
the choice of $L$ and $\eps$, we have \be \nonumber \frac{\log
\beta_1}{\log \beta_*} \geq 1 - \frac{1}{2N}, \ee thus $\rho_*$
has H\"older-exponent of continuity $1-\frac{1}{2N}$. Define \be
\nonumber G_N^* := \{\uz \in \Sigma^*: \exists \nu \in
V_\sigma(\uz) \text{ s.t. } h(\nu) < (1-1/N) \log\beta_*\}. \ee As
before, we have $\Omega_N \cap [\beta_1,\beta_2] \subset
\rho_*(G_N^* \cap D^*)$ and $\htop(G_N^*) \leq
(1-1/N)\log\beta_*$. Finally $\dim_H (\Omega_N \cap
[\beta_1,\beta_2]) < 1$ and $\lambda(\Omega_N \cap
[\beta_1,\beta_2]) = 0$. $\Box$

\textbf{Proof of the Corollary:} The proof is similar to the proof
of Corollary \ref{cl selfnormality for Tab alpha fixed}. Equation
\eqref{eq T^n = phi sigma^n i} is true, since we work on $B_0$.
$\Box$

In particular, when we consider the tent map ($s = (1,-1)$), we recover the main Theorem of
Bruin in \cite{Br98}. We do not state this theorem for all $x \in [0,1]$ as for the map
$\Tab$, because we do not have an equivalent of Lemma \ref{lm separation of coding in Tb}
for all $x \in [0,1]$. This is the unique missing step of the proof.

\noindent
{\bf Acknowledgements:\enspace} We thank H. Bruin for
correspondence about Proposition 1 and for communicating us
results before publication.

\end{document}